\documentclass{proc-l-arxiv}
\usepackage{amsmath, amsthm, amscd, amsfonts, amssymb, graphicx, color}
\newtheorem{theorem}{Theorem}[section]
\newtheorem{lemma}[theorem]{Lemma}
\newtheorem{corollary}[theorem]{Corollary}

\theoremstyle{definition}
\newtheorem{definition}[theorem]{Definition}

\theoremstyle{remark}
\newtheorem{remark}[theorem]{Remark}

\numberwithin{equation}{section}

\def\be{\begin{equation}}
\def\ee{\end{equation}}

\newcounter{alphabet}
\newcounter{tmp}

\makeatletter
\newcommand{\Ref}[1]{\@ifundefined{r@#1}{}{\setcounter{tmp}{\ref{#1}}\Alph{tmp}}}
\makeatother
\newcommand{\IN}{{\mathbb N}}
\newcommand{\IC}{{\mathbb C}}
\newcommand{\ID}{{\mathbb D}}

\newcommand{\bee}{\begin{enumerate}}
\newcommand{\eee}{\end{enumerate}}
\newcommand{\bthm}{\begin{theorem}}
\newcommand{\ethm}{\end{theorem}}
\newcommand{\blem}{\begin{lemma}}
\newcommand{\elem}{\end{lemma}}
\newcommand{\bdefe}{\begin{definition}}
\newcommand{\edefe}{\end{definition}}
\newcommand{\bcor}{\begin{corollary}}
\newcommand{\ecor}{\end{corollary}}
\newcommand{\br}{\begin{remark}}
\newcommand{\er}{\end{remark}}
\newcommand{\bpf}{\begin{proof}}
\newcommand{\epf}{\end{proof}}

\newcommand{\ba}{\begin{array}}
\newcommand{\ea}{\end{array}}
\newcommand{\beq}{\begin{eqnarray}}
\newcommand{\beqq}{\begin{eqnarray*}}
\newcommand{\eeq}{\end{eqnarray}}
\newcommand{\eeqq}{\end{eqnarray*}}

\newcommand{\ds}{\displaystyle}
\begin{document}

\title[On the Bohr inequality with a fixed zero coefficient]{On the Bohr inequality with a fixed zero coefficient}


%    author one information

\author[S.A. Alkhaleefah]{Seraj A.~Alkhaleefah}

\address{Seraj A.~Alkhaleefah, Kazan  Federal University, 420 008 Kazan, Russia
}
\email{s.alkhaleefah@gmail.com}

\author[I. R Kayumov]{Ilgiz R Kayumov}
\address{Ilgiz R Kayumov, Kazan  Federal University, 420 008 Kazan, Russia
}
\email{ikayumov@gmail.com}

\author[S. Ponnusamy]{Saminathan Ponnusamy
%$^\dagger $
%${}^{~\mathbf{*}}$
}
\address{Saminathan  Ponnusamy,
Department of Mathematics,
Indian Institute of Technology Madras,
Chennai-600 036, India.}
\email{samy@iitm.ac.in}

\thanks{}

%    \subjclass is required.


\date{}

\dedicatory{}

%    "Communicated by" -- provide editor's name; required.
\commby{}

%    Abstract is required.
\begin{abstract}
In this paper, we introduce the study of the Bohr phenomenon for a quasi-subordination family of functions, and
establish the classical Bohr's inequality for the class of quasisubordinate functions.
As a consequence, we improve and obtain the exact version of the classical Bohr's inequality for bounded analytic functions and also for
$K$-quasiconformal harmonic mappings by replacing the constant term by the absolute value of the analytic part of the
given function. We also obtain the Bohr radius for the subordination family of odd analytic functions.
\end{abstract}

\maketitle

%%%%%%%%%%%%%%%%%%%%%%%%%%%%%%%%%%%%%%%%%%%%%%%%%%%%%%%%%%%%%%%%%%%%%%%%%%%%%%%%%%%%%%%%%%%
\section{Introduction and Preliminaries}

In this article, our primary concern is to study Bohr's phenomenon for the class of quasi-subordination functions and obtain the exact version of the classical Bohr's inequality for the case of analytic functions
and also for the case of harmonic functions defined on the open unit disk $\ID=\{z\in \IC:\, |z|<1\}$.
The classical result of H. Bohr \cite{Bohr-14}, which in the final form was proved
independently by M. Riesz, I. Schur and N. Wiener, is as follows:

%\begin{Thm}\label{BohrTh1}

\vspace{8pt}
\noindent
{\bf Theorem A.}
{\it
Let $f(z)=\sum_{k = 0}^\infty a_kz^k$ be analytic in $\ID$ and $ |f(z)| \le 1 $ for all $z \in \ID$. Then
\begin{equation}\label{bohr1}
|a_0| +\sum_{k=1}^\infty |a_k|r^k \le 1\  \ \text{for all}\  \ r \le \frac{1}{3}
\end{equation}
and the constant $1/3$, called the Bohr radius, cannot be improved.
}
%\end{Thm}
\vspace{8pt}

In 1956, Ricci  \cite{Ricci-1956} initiated the investigation of the Bohr radius with fixed zero-coefficient $a_0$,
and in 1962, Bombieri \cite{Bom-62} solved the problem for $|a_0|\ge \frac12$.  In this paper in the later part of our
investigation (see Theorem \ref{AKP-th5}), we strengthen these
results and furthermore, in Theorem \ref{AKP-th6}, we extend it for sense-preserving $K$-quasiconformal harmonic mappings of the unit disk.

In the recent years, the problem about the Bohr radius attracted the attention of many researchers in various directions in
functions of one and several complex variables:
to planar harmonic mappings, to polynomials, to domains in several complex variables,  to solutions of elliptic partial
differential equations and to more abstract settings. For more information about Bohr's inequality stated above and further
related investigations, we refer the reader to the recent survey articles on the Bohr radius from
\cite{AAPon1, IsmaKayKayPon1}, \cite[Chapter 8]{GarMasRoss-2018},
%the articles \cite{PaulPopeSingh-02-10,PaulSingh-04-11,Sidon-27-15,Tomic-62-16}
and the references therein. See also \cite{Bom-62,BomBor-04}. In particular,
Boas and Khavinson \cite{BoasKhavin-97-4},  Aizenberg \cite{Aizen-00-1,Aiz07}, and Aizenberg and  Tarkhanov \cite{bib:10} have extended the
Bohr inequality  for holomorphic functions on certain specific domains (such as complete Reinhardt domain) in $\IC^n$.

Recently,  Kayumov et al. \cite{KayPonShak1} investigated Bohr's radius for   locally univalent planar harmonic mappings.
Several improved versions of the classical Bohr's inequality were given by Kayumov and Ponnusamy in  \cite{KayPon3}  (see also \cite{KayPon2})
whereas Evdoridis et al.  \cite{EvPoRa-2017} have presented several improved versions of Bohr's inequality for harmonic mappings.
In \cite{KayPon1}, Kayumov and Ponnusamy also discussed Bohr's radius for the class of analytic functions $g$, when $g$ is subordinate to a member of the class of odd univalent functions. For certain recent results, we refer to \cite{AliBarSoly,BDK5,KayPon2}.
%\cite{Abu2,Abu4,Abu3,AliBarSoly,BDK5,KayPon2}.
In particular, Kayumov and Ponnusamy \cite{KayPon2} established
the following theorem which settled the open problem proposed by Ali et al. \cite{AliBarSoly}.

%\begin{Thm}\label{KayPonn}

\vspace{8pt}
\noindent
{\bf Theorem B.}
{\it
If a function $f(z)=\sum_{k = 1}^\infty a_{2k-1}z^{2k-1}$ is odd analytic in $\ID$ and $|f(z)|\le 1$ in $\ID$, then
$$%\begin{equation}\label{eq-KayPonn}
\sum_{k=1}^\infty |a_{2k-1}|r^{2k-1} \le 1 \  \ \text{for all}\  \ r \le r_0,
$$%\end{equation}
where $r_0\simeq0.789991...$ is the maximal positive root of the equation
$$ 8r^4 +r^2 - 6r +1 = 0 $$
and the constant $r_0$ cannot be improved.
}
%\end{Thm}
\vspace{8pt}

In 1970, Robertson \cite{Robert-70} introduced and developed
the concept of quasi-subordination which combines the principles of subordination and majorization.

If  $f$ and $g$ are analytic in $\ID$, $\omega$ is a Schwarz function (i.e.  $\omega$ is analytic in $\ID$,
$\omega(0)=0$ and $|\omega (z)|\leq 1$ for $|z|<1$) and all three satisfy $f(z) = g(\omega(z))$ for $z\in \ID$, then
we write $f(z)\prec g(z)$ in $\ID$ and say that $f$ is \textit{subordinate} to $g$. The importance of the principle of subordination
stems from the fact that when $f$ is subordinate to $g$, $f(\ID) \subset g(\ID)$ and this has been extensively used in the literature.
We say that $f(z)$ is \textit{majorized} by $g(z)$ in $\ID$ if $|f(z)|\leq |g(z)|$ for all $z\in\ID$.
%, and we denote this by writing $f(z) \ll g(z)$.

\bdefe
For any two analytic functions $f$ and $g$ in $\ID$, we say that the function $f$ is \textit{quasi-subordinate} to $g$ (relative to $\Phi$), denoted by $f(z)\prec _q g(z)$ in $\ID$, if
there exist two functions $\Phi $  and  $\omega$, analytic in $\ID$, satisfying $\omega(0)=0$, $|\Phi(z)|\leq 1$ and $|\omega (z)|\leq 1$ for $|z|<1$ such that
\begin{equation}\label{Eq7a}
f(z)=\Phi(z)g(\omega(z)).
\end{equation}
\edefe

There are two special cases which are of particular interest. The choice $\Phi(z)=1$ corresponds to subordination, whereas
$\omega(z) = z$ gives majorization, i.e. \eqref{Eq7a} reduces to the form $f(z)=\Phi (z)g(z)$. In other words, if either $f\prec g$ or $|f(z)|\leq |g(z)|$ in $\ID$, then $f(z)\prec _q g(z)$ in $\ID$.
Thus, the notion of quasi-subordination generalizes both the concept of subordination and the principle of majorization. Several theorems exist in the literature that relate with these two concepts and are widely used in function theory, and some of the known results continue to hold in the setting of quasi-subordination. See \cite{MacGre-67, Rogo-43}. Note also that \eqref{Eq7a} is equivalent to saying that the quotient $f(z)/\Phi(z)$ is analytic and is subordinate to $g(z)$ in $\ID$.

\br
{\it
 On the linear space ${\mathcal H}(\ID)$ consisting of complex-valued analytic functions $g$  defined on $\ID$, let $\omega$
denote self-map of $\ID$.
The composition operator $C_\omega$ with symbol $\omega$ is defined as
$$ C_\omega g = g\circ\omega ~ \mbox{ for }~  g \in {\mathcal H}(\ID).
$$
Similarly, a weighted composition operator $W_{\omega, \Phi}$ is an operator that maps $g\in{\mathcal H}(\ID)$ into
$W_{\omega, \Phi}(f)=\Phi(z)g(\omega(z))$, where $\Phi $  and $\omega$ are analytic defined on $\ID$ such that $\omega (\ID) \subset \ID$.
Note also that for a given complex-valued function $\Phi$ defined on $\ID$, the multiplication operator with symbol $\omega$ is defined by
$$M_\Phi g = \Phi g  ~ \mbox{ for }~  g \in {\mathcal H}(\ID).
$$
These operators appear in a natural way, for example,  in the study of a number of questions about the boundedness and compactness of operators on
various function spaces in a more general setting. See \cite{CowenMach95}.
Thus, it is worth pointing out that the concept of \emph{`subordination'} is nothing but a composition (operator) with a function mapping $\ID$ into itself,
and the concept of\emph{`quasi-subordination'}
is nothing but a weighted composition (operator). Note also that the multiplication operator is related to the majorization.
}
\er

The paper is organized as follows.
Section \ref{AKP-sec5} is devoted to state our main results whose proofs will be presented in Section \ref{AKP-sec2}.
First we show that  (Theorem \ref{AKP-th4}) the radius $\frac13$ of Bohr inequality remains the same even when the functions $f$ and $g$ are related with a quasi-subordination relation \eqref{Eq7a} which clearly reveals the fact that the classical Bohr inequality continues to hold in a more general setting. Secondly, as a consequence of Theorem \ref{AKP-th4}, we  present in  Corollary \ref{AKP-th2} the exact version of Theorem~A. % \Ref{BohrTh1}.
Thirdly, we show in Theorem \ref{AKP-th1} that the Bohr radius for the subordinating family of odd functions is $1/{\sqrt3}$. In Theorem \ref{AKP-th3}, we present a sharp version of Bohr's inequality for sense-preserving $K$-quasiconformal harmonic mappings. Finally, in Theorems \ref{AKP-th5} and \ref{AKP-th6}, we essentially investigate the Bohr phenomenon by replacing the constant term by the function  itself in the case of analytic functions, and by the analytic part  in the case of harmonic functions, respectively.

%Bohr originally discovered this inequality \eqref{bohr1}, known as the Bohr inequality, only for $r\le 1/6$. The number  $1/3$ is called the Bohr radius for the class of bounded analytic functions in $\ID$.
%%But later the fact that the inequality is actually true for $r\le 1/3$, and the value $1/3$ is sharp, was obtained independently by M. Riesz, I. Schur and N. Weiner.
%Moreover, equality in $\sum_{k=0}^\infty |a_k|(1/3)^k \le 1$ holds for constant functions only.

\section{Main Results and their consequences}\label{AKP-sec5} %Bohr's inequality for quasi-subordinating family of functions
First we state an improved version of Bohr's inequality for a quasi-subordinating family of functions.

\bthm\label{AKP-th4}
Let $f(z)$ and $g(z)$ be two analytic functions in $\ID$ with the Taylor series expansions
$f(z)=\sum_{k = 0}^\infty a_kz^k$ and $g(z)=\sum_{k = 0}^\infty b_kz^k$ for $z\in\ID$. If $f(z)\prec_q g(z)$,  then
$$\sum_{k=0}^\infty |a_k| r^k \le \sum_{k=0}^\infty |b_k| r^k \  \ \text{for all}\  \ r \le \frac{1}{3}.
$$
%for $|z|=r\le \frac{1}{3}$.
\ethm

We are now ready to state two simple corollaries which are of independent interest, and the first of which was obtained
recently by  Bhowmik and Das \cite{BhowDas-18}.
%Before we continue the discussion, let us state the corollaries  which also play a crucial role in establishing
%improved and extended versions of Theorem \Ref{BohrTh1}.

\bcor\label{Bhowmik}
{\rm \cite{BhowDas-18}}
Let $f(z)$ and $g(z)$ be two analytic functions in $\ID$ such that  $f(z)=\sum_{k = 0}^\infty a_kz^k$, and
$g(z)=\sum_{k = 0}^\infty b_kz^k$. If $f(z)\prec g(z)$ in $\ID$, then
$$%\begin{equation}\label{eq-Bhowmik}
\sum_{k=0}^\infty |a_k|r^k \le \sum_{k=0}^\infty |b_k|r^k \  \ \text{for all}\  \ r \le \frac{1}{3}
$$%\end{equation}
and the constant $1/3$ cannot be improved.
\ecor

\bcor\label{Analog_Bhowmik}
Let $f(z)=\sum_{k = 0}^\infty a_kz^k$, and $g(z)=\sum_{k = 0}^\infty b_kz^k$  be two analytic functions in $\ID$. If
$f(z)$ is majorized by $g(z)$, i.e. $|f(z)| \leq |g(z)|$ in $\ID$, then
$$
\sum_{k=0}^\infty |a_k|r^k \le \sum_{k=0}^\infty |b_k|r^k \  \ \text{for all}\  \ r \le \frac{1}{3}
$$
and the constant $1/3$ cannot be improved.
\ecor

These corollaries play a crucial role in establishing
generalized versions of Theorem~A. %\Ref{BohrTh1}.

According to Theorem~B, %\Ref{KayPonn},
the Bohr radius for the class of odd functions is $0.789991...$ and thus, it is natural to ask for the Bohr radius for the subordinating family of odd analytic functions. Our next result answers this question.

\bthm\label{AKP-th1}
Let $f(z)$ and $g(z)$ be odd analytic functions in $\ID$ with Taylor expansions
$f(z)=\sum_{k = 1}^\infty a_{2k-1} z^{2k-1}$ and $g(z)=\sum_{k = 1}^\infty b_{2k-1} z^{2k-1},$
%$$
%f(z)=\sum_{k = 1}^\infty a_{2k-1} z^{2k-1} ~\mbox{ and }~   g(z)=\sum_{k = 1}^\infty b_{2k-1} z^{2k-1},
%$$
respectively. If $f(z)\prec g(z)$ , then
\begin{equation}\label{AKP-eq6}
\sum_{k=1}^\infty |a_{2k-1}| r^{2k-1} \le \sum_{k=1}^\infty |b_{2k-1} | r^{2k-1}
~\mbox{ for $\ds |z|=r\le \frac1{\sqrt3}$}.
\end{equation}
%for $|z|=r\le \frac1{\sqrt3}$.
\ethm

A harmonic mapping in $\ID$ is a complex-valued function $f$ in $\ID$, which satisfies the Laplace equation
$\Delta f=4f_{z\,\overline{z}}=0$.  It follows that $f$ admits the canonical representation
$f= h+\overline{g}$, where $h$ and $g$ are analytic in $\ID$ with $f(0)=h(0)$.
The Jacobian $J_{f}$ of $f$ is given by $J_{f} =|h'|^2-|g'|^2.$ We say that $f$ is sense-preserving in $\ID$ if $J_{f}(z)>0$ in $\ID$. Consequently,
$f$ is locally univalent and sense-preserving in $\ID$ if and only if
$J_{f}(z)>0$ in $\ID$; or equivalently if $h'\neq 0$ in $\ID$ and the dilatation
$\omega_f=:\omega =g'/h'$ has the property that $|\omega (z)|<1$ in $\ID$. For a detailed treatment of the geometric point of view of
planar harmonic mappings of the unit disk, we refer to \cite{Duren:Harmonic} and also \cite{Clunie-Small-84,SaRa2013}.

In order to state  our result about the Bohr radius for quasiconformal harmonic mappings, we recall that
a sense-preserving homeomorphism $f$ from the unit disk $\ID$ onto $\Omega'$, contained in the Sobolev class $W_{loc}^{1,2}(\ID)$,
is said to be a {\it $K$-quasiconformal mapping} if, for $z\in\ID$,
$$\frac{|f_{z}|+|f_{\overline{z}}|}{|f_{z}|-|f_{\overline{z}}|}= \frac{1+|\omega_f(z)|}{1-|\omega_f(z)|}\leq K,
 ~\mbox{ i.e., }~|\omega_f(z)|=\left |\frac{g'(z)}{h'(z)}\right |\leq k=\frac{K-1}{K+1},
$$
where $K\geq1$ so that $k\in [0,1)$.
We now state a new version of Bohr's inequality for harmonic mappings.

\bthm\label{AKP-th3}
Suppose that $f(z)= h(z) + \overline{g(z)} =\sum_{n=0}^\infty a_nz^n + \overline{\sum_{n = 1}^\infty b_nz^n}$ is a sense-preserving
$K$-quasiconformal harmonic mapping of the disk $\ID$ (or more generally $|\omega_f(z)| \leq k$ for some $k\in [0,1]$), where $|h(z)|\le 1$  in $\ID$. Then the following sharp inequalities hold:
\begin{equation}\label{Eq6}
\frac{1-r(|a_0|+(k+1)(1-|a_0|^2))}{1-r|a_0|}+\sum_{n=1}^\infty |a_n|r^n + \sum_{n=1}^\infty |b_n|r^n  \le 1
~\mbox{ for  $\ds r \le \frac{1}{3}$}.
%~\mbox{ for  $\ds r \le \frac{K+1}{5K+1}$}.
\end{equation}
 The functions
$$
\frac{z+a_0}{1+\overline{a_0}z}+\lambda \overline{\frac{z+a_0}{1+\overline{a_0}z}}
$$ with $\lambda \to 1$ demonstrate that the inequality \eqref{Eq6} is sharp for all $a_0 \in \ID$ and all
$r \leq \frac{1}{3}$.
%$r \leq \frac{K+1}{5K+1}$.
\ethm

Before we continue the discussion, let us remark that the classical Bohr inequality is not sharp for any
individual function. Namely, it is easy to show that for any given function
the Bohr radius is always greater than  $1/3$. As a result of Theorem \ref{AKP-th3}, here is the sharp result
which shows that $1/3$ cannot be improved even in the case of individual functions.

\bcor\label{AKP-th2}
Let $f(z)=\sum_{n=0}^\infty a_nz^n$ be an analytic function in $\ID$ and $|f(z)|\leq 1$ for all $z \in \ID$. Then the following
sharp inequality holds:
\begin{equation}\label{AKP-eq7}
\frac{1-(1+|a_0|-|a_0|^2)r}{1-|a_0|r} +\sum_{n=1}^\infty |a_n|r^n \le 1 ~\mbox{ for all $r \le 1/3$},
\end{equation}
The function $g(z)=(z+a_0)/(1+\overline{a_0}z)$ shows that equality holds for all $a_0 \in \ID$ and $r \leq 1/3$ .
\ecor
\bpf
The result follows if we let $K=1$ (i.e.  $k=0$) in Theorem \ref{AKP-th3} so that $g(z)\equiv 0$ in $\ID$.
%$K\rightarrow 1^{+}$ (i.e.  $k\rightarrow 0^{+}$) in Theorem \ref{AKP-th3}.
\epf

%\bcor
%Suppose that $f(z)= h(z) + \overline{g(z)} =\sum_{n=0}^\infty a_nz^n + \overline{\sum_{n = 1}^\infty b_nz^n}$ is a sense-preserving
%$K$-quasiconformal harmonic mapping of the disk $\ID$, where $|h(z)|\le 1$  in $\ID$. Then the following sharp inequalities hold:
%\begin{equation}\label{Eq6}
%\frac{1-r(|a_0|+k(1-|a_0|^2))}{1-r|a_0|}+\sum_{n=1}^\infty |a_n|r^n + \sum_{n=1}^\infty |b_n|r^n  \le 1
%~\mbox{ for  $\ds r \le \frac{K+1}{5K+1}$}.
%\end{equation}
% The functions
%$$
%\frac{z+a_0}{1+\overline{a_0}z}+\lambda \overline{\frac{z+a_0}{1+\overline{a_0}z}}
%$$ with $\lambda \to 1$ demonstrate that the inequality \eqref{Eq6} is sharp for all $a_0 \in \ID$ and all $r \leq \frac{K+1}{5K+1}$.
%\ecor

\bthm\label{AKP-th5}
Suppose that $f(z)=\sum_{k = 0}^\infty a_kz^k$ is an analytic function in $\ID$ and $|f(z)|<1$ for all $z \in \ID$ and $0\le |a_0|=a < 1$. Then
\begin{equation}\label{AKP-eq8}
|f(z)|+\sum_{k=1}^\infty |a_k|r^k \le 1
\end{equation}
for all $a\ge 2\sqrt3-3 \approx 0.4641016$ and $|z|=r\le r_a$, where
$$
r_a = \frac{\sqrt{(1+a)^2+a^2}-(1+a)}{a^2} =\frac{1}{\sqrt{(1+a)^2+a^2}\, +1+a}
$$
and the radius $r_a$ is sharp.
\ethm

\br \label{AKP-rem1}
From the proof of Theorem \ref{AKP-th5}, it can be easily seen that for $r\le \sqrt5-2$ the inequality \eqref{AKP-eq8}
continues to hold for all $a< 1$.
\er

We now generalize Theorem \ref{AKP-th5} in order to present a generalized version of Bohr's inequality with constant term by analytic
part of the corresponding harmonic function.

\bthm\label{AKP-th6}
Suppose that $f(z)= h(z) + \overline{g(z)} =\sum_{n = 0}^\infty a_nz^n + \overline{\sum_{n = 1}^\infty b_nz^n}$ is a sense-preserving
$K$-quasiconformal harmonic mapping of the disk $\ID$, where $|h(z)|< 1$  in $\ID$ and $0\le a=|a_0| < 1$.
Then the following sharp inequalities hold:
\begin{equation}\label{Eq9}
|h(z)|+\sum_{n=1}^\infty |a_n|r^n + \sum_{n=1}^\infty |b_n|r^n  \le 1
\end{equation}
for all $a\ge \alpha _k$ and $|z|=r\le r_{a,k}$, where
$$ \alpha _k = \frac{\sqrt{k^2+12k+12}-(2k+3)}{k+1}
~\mbox{ and }~r_{a,k} =\frac{B_{a,k}-(k+2)(1+a)}{2a^2(k+1)+2ak}
$$
with
$B_{a,k}= \sqrt{a^2(k^2+8k+8)+2a(k^2+6k+4)+(k+2)^2}.
$
The radius $r_{a,k}$ is sharp.
\ethm

 \section{Proofs of the Main Results}\label{AKP-sec2}
In the proofs of Theorems \ref{AKP-th4} and   \ref{AKP-th1} , we will
use some approaches used in \cite{Rogo-43} (see also
  \cite[Proof of Lemma 1]{BhowDas-18}) .

\subsection{Proof of Theorem \ref{AKP-th4}}
Suppose that  $f\prec_q g$. Then there exist two analytic functions $\Phi$ and $\omega$ satisfying $\omega(0)=0$, $|\omega(z)|\le1$ and
$|\Phi(z)|\le1$ for all $z\in \ID$ such that %\eqref{Eq7a} holds
\begin{equation}\label{Eq7}
f(z)=\Phi(z)g(\omega(z)).
\end{equation}
Now for the analytic function $\omega(z)=\sum_{n=1}^\infty \alpha_n z^n$, the Taylor expansion of the
$k$-th power of $\omega$, where $k \in \IN$, can be written as
\begin{equation}\label{Eq8}
\omega^k(z) = \sum_{n=k}^\infty \alpha_n^{(k)} z^n.
\end{equation}
We observe that, since $\omega(0)=0$ and $|\omega(z)|\le1$, it follows from Theorem~A %\Ref{BohrTh1}
that
\be\label{EQ8b}
\sum_{n=k}^\infty|\alpha_n^{(k)}|r^{n-k}\le1 \ \ \text{for all}\ \ r\le\frac13 .
\ee
For the analytic function $\Phi(z)$, we may write $\Phi(z)=\sum_{m=0}^\infty \phi_m z^m$ and thus, by Theorem~A, %\Ref{BohrTh1},
we have
\be\label{EQ8a}
\sum_{m=0}^\infty |\phi_m|r^m\le1 \ \ \text{for all}\ \ r\le\frac13 .
\ee
Also, from the equality \eqref{Eq7}, taking into consideration from \eqref{Eq8} that
$$
\omega^0(z)=1=\sum_{n=0}^\infty \alpha_n^{(0)}z^n, \ \ \text{where}\ \ \ \alpha_0^{(0)}=1, ~ \alpha_n^{(0)}= 0 \mbox{ for $n\geq 1$},
$$
we can rewrite the quasi-subordinate relation \eqref{Eq7a} with the help of \eqref{Eq8} in series form as
\beqq
\sum_{k=0}^\infty a_k z^k %&=& \sum_{m=0}^\infty \phi_m z^m \sum_{k=0}^\infty b_k \omega^k(z)\\
&=& \sum_{m=0}^\infty \phi_m z^m \left(\sum_{k=0}^\infty b_k \sum_{n=k}^\infty \alpha _n^{(k)} z^n\right)\\
&=& \sum_{m=0}^\infty \phi_m z^m \left(\sum_{k=0}^\infty \left (\sum_{n=0}^k b_n\alpha_k^{(n)}\right ) z^k\right)\\
&=& \sum_{m=0}^\infty \phi_m z^m \sum_{k=0}^\infty B_k z^k,% \ \ \text{where}\ \ \ B_k = \sum_{n=0}^k b_n\alpha_k^{(n)},\\
%&=& \sum_{k=0}^\infty \left(\sum_{m+j=k}\phi_m B_j\right) z^k.
\eeqq
where $B_k = \sum_{n=0}^k b_n\alpha_k^{(n)}$.
%, and so,
%\be\label{EQ8c}
%|B_k| \leq  \sum_{n=0}^k |b_n|\, |\alpha_k^{(n)}|
%\ee
Thus, the last relation takes the form
$$\sum_{k=0}^\infty a_k z^k = \sum_{k=0}^\infty \left(\sum_{m+j=k}\phi_m B_j\right) z^k,
$$
which by equating the coefficients of $z^k$ on both sides gives % we have, for any $k \ge 0$:
\begin{equation}\label{Eq10}
a_k = \sum_{m+j=k}\phi_m B_j ~\mbox{ for each $k \ge 0$}.
%\ \ \text{where}\ \ \ B_j = \sum_{n=0}^j b_n\alpha_j^{(n)}
\end{equation}
Applying the triangle inequality to the last relation shows that
\beqq
\sum_{k=0}^\infty |a_k|r^k %&=& \sum_{k=0}^\infty \left|\sum_{m+j=k}\phi_m B_j\right| r^k \\
&\le & \sum_{k=0}^\infty \left (\sum_{m+j=k}|\phi_m| \,|B_j| \right ) r^k
%\\&=&
=\sum_{k=0}^\infty \sum_{m+j=k}|\phi_m|r^m |B_j|r^j \\
&= &
\left (\sum_{m=0}^\infty |\phi_m|r^m \right )\sum_{k=0}^\infty|B_k|r^k\\
&\le& \sum_{k=0}^\infty|B_k|r^k ~\mbox{ for all $r\le\frac13$, (by \eqref{EQ8a})}.%\ \ \left(\text{where}\ \ \ \sum_{m=0}^\infty |\phi_m|r^m\le1 \ \ \text{for all}\ \ r\le\frac13\right)\\
\eeqq
Also, because $|B_k| \leq  \sum_{n=0}^k |b_n|\, |\alpha_k^{(n)}|$, we obtain that
\beqq
\sum_{k=0}^\infty|B_k|r^k %&=& \sum_{k=0}^\infty\left|\sum_{n=0}^k b_n\alpha_k^{(n)}\right|r^k \\
&\le & \sum_{k=0}^\infty\sum_{n=0}^k |b_n|\, |\alpha_k^{(n)}|r^k
%\\&=&
=\sum_{k=0}^\infty|b_k| \sum_{n=k}^\infty|\alpha_n^{(k)}|r^n \\
&= & \sum_{k=0}^\infty|b_k| \left (\sum_{n=k}^\infty|\alpha_n^{(k)}|r^{n-k}\right ) r^k\\
&\le& \sum_{k=0}^\infty|b_k|r^k  ~\mbox{ for all $r\le\frac13$, (by \eqref{EQ8b})}%\ \ \left(\text{where}\ \ \ \sum_{n=k}^\infty|\alpha_n^{(k)}|r^{n-k}\le1 \ \ \text{for all}\ \ r\le\frac13\right)
\eeqq
and hence, we obtain that
$$
\sum_{k=0}^\infty |a_k|r^k \le \sum_{k=0}^\infty|B_k|r^k \le \sum_{k=0}^\infty|b_k|r^k \ \ \text{for all}\ \ r\le\frac13.
$$
The proof of Theorem \ref{AKP-th4} is complete. $\hfill \Box$

\subsection{Improved version of the classical Bohr inequality for odd and $p$-symmetric functions}

%We recall that the radius $\frac13$ of Bohr inequality remains the same even when the functions $f$ and $g$ are related with
%quasi-subordination relation  (Theorem \ref{AKP-th4}) and also for the subordinating family of functions (Corollary \ref{Bhowmik}).
%Compare with the work of Kayumov and Ponnusamy in Theorem \Ref{KayPonn} for the class of odd functions where  the Bohr
%radius turned out to be $0.789991...$. Our next aim is to show that the Bohr radius for the odd functions
%$f$ and $g$ with subordination relation turned out to be $\frac1{\sqrt3}$.

Our next result is indeed a simple consequence of Corollary \ref{Bhowmik} and we state it in this form because of its
independent interest.

\blem\label{AKP-lem1}
Let $f(z)$ and $g(z)$ be analytic and $p$-symmetric in $\ID$  with the Taylor expansions
$f(z)=\sum_{k = 0}^\infty a_{pk} z^{pk} $  and $ g(z)=\sum_{k = 0}^\infty b_{pk} z^{pk},$
%$$
%f(z)=\sum_{k = 0}^\infty a_{pk} z^{pk} ~\mbox{ and }~  g(z)=\sum_{k = 0}^\infty b_{pk} z^{pk},
%$$
respectively. If $f(z)\prec g (z)$, then
$$\sum_{k=0}^\infty |a_{pk}| r^{pk} \le \sum_{k=0}^\infty |b_{pk} | r^{pk} \ \ \text{for all}\ \ r\le \frac1{\sqrt[p]3}.
$$
The constant  $1/{\sqrt[p]3}$ cannot be improved.
\elem
\bpf
It suffices to set $\zeta =z^p$, and consider the functions
$$
f_1(\zeta)= \sum_{k=0}^\infty a_{pk}\zeta^k ~\mbox{ and }~ g_1(\zeta)= \sum_{k=0}^\infty b_{pk}\zeta^k.
$$
Then $f_1(\zeta)\prec g_1(\zeta)$ for $|\zeta|<1$ and, by Corollary \ref{Bhowmik},  we obtain that
$$
\sum_{k=0}^\infty |a_{pk}|\,|\zeta|^k \le \sum_{k=0}^\infty |b_{pk}|\,|\zeta|^k \ \ \text{for all} \ \ |\zeta| =|z|^p\le \frac13.
$$
The desired conclusion follows.
\epf

\subsection{Proof of Theorem \ref{AKP-th1}}
 Let $f\prec g$, where $f$ and $g$ are as in the statement. Then there exists a function $\omega$, analytic in $\ID$, satisfying $\omega(0)=0$ and $|\omega(z)|\le1$ for all $|z|<1$ such that $f(z)=g(\omega(z))$ which in terms of series can be written as
\begin{equation}\label{AKP-eq3}
\sum_{k=1}^\infty a_{2k-1} z^{2k-1} = \sum_{k=1}^\infty b_{2k-1} \omega(z)^{2k-1},
\end{equation}
where, as usual, we write $\omega(z)=\sum_{n=1}^\infty \alpha_n z^n$ and the Taylor expansion of the $(2k-1)$-th power of $\omega$,
where $k \in \IN$, has the form
%$$
%\omega^k(z) = \sum_{n=k}^\infty \alpha_n^{(k)} z^n
%$$
%So
\begin{equation}\label{AKP-eq4}
\omega^{2k-1}(z) = \sum_{n=2k-1}^\infty \alpha_n^{(2k-1)} z^n.
\end{equation}
Now we plug the equality \eqref{AKP-eq4} into the right hand side of the relation \eqref{AKP-eq3}, and obtain
\beqq
\sum_{k=1}^\infty a_{2k-1}z^{2k-1} &=& \sum_{k=1}^\infty b_{2k-1}\left (\sum_{n=2k-1}^\infty \alpha_n^{(2k-1)} z^n \right ) \\
%&=&  b_1 \sum_{n=1}^\infty \alpha_n^{(1)} z^n + b_3 \sum_{n=3}^\infty \alpha_n^{(3)} z^n
%+ b_5 \sum_{n=5}^\infty \alpha_n^{(5)} z^n + \cdots\\
&=& \left ( b_1\alpha_1^{(1)}z + b_1\alpha_2^{(1)}z^2 + b_1\alpha_3^{(1)}z^3 + \cdots \right ) + \left ( b_3\alpha_3^{(3)}z^3 +  b_3\alpha_4^{(3)}z^4 +\right .\\
 && \left . + b_3\alpha_5^{(3)}z^5 + \cdots \right ) + \left ( b_5\alpha_5^{(5)}z^5 + b_5\alpha_6^{(5)}z^6 + b_5\alpha_7^{(5)}z^7 + \cdots \right ) + \cdots .
\eeqq
Clearly, the coefficients of $z^{2n}$ have to be zero and thus, $\ds \alpha_{2m}^{(2k-1)} =0$ for $m=1,2, \ldots$.
%$$\alpha_2^{(2k-1)}=\alpha_4^{(2k-1)}= \cdots  =\alpha_{2n}^{(2k-1)} = \cdots  = 0
%$$
Thus, we can write the last equation as
$$\sum_{k=1}^\infty a_{2k-1}z^{2k-1} = \sum_{k=1}^\infty \left (\sum_{n=1}^k b_{2n-1}\alpha_{2k-1}^{(2n-1)} \right ) z^{2k-1}
$$
and equating the coefficients of $z^{2k-1}$ on both sides, we have
$$a_{2k-1} = \sum_{n=1}^k b_{2n-1} \alpha_{2k-1}^{(2n-1)} ~\mbox{ for any $k \ge 1$}.
$$
Applying the triangle inequality to the last relation shows that
%$$
%|a_{2k-1}| = \sum_{n=1}^k | b_{2n-1}|\, \left | \alpha_{2k-1}^{(2n-1)}\right | ~\mbox{ for any $k \ge 1$}.
%$$
\beq
\sum_{k=1}^m|a_{2k-1}|r^{2k-1}
%&=& \sum_{k=1}^m |\sum_{n=1}^k b_{2n-1}\alpha_{2k-1}^{2n-1}|r^{2k-1}\\
&\le& \sum_{k=1}^m \left (\sum_{n=1}^k |b_{2n-1}|\,\left |\alpha_{2k-1}^{2n-1}\right | \right )r^{2k-1} \label{AKP-eq5} \\
%&=& |b_1||\alpha_1^{(1)}|r + \sum_{n=1}^2|b_{2n-1}||\alpha_3^{(2n-1)}|r^3 + \cdots +\sum_{n=1}^m|b_{2n-1}||\alpha_m^{(2n-1)}|r^m\\
%&=& |b_1||\alpha_1^{(1)}|r + |b_1||\alpha_3^{(1)}|r^3 + |b_3||\alpha_3^{(3)}|r^3 + \cdots  + |b_1||\alpha_m^{(1)}|r^m + \cdots + |b_m||\alpha_m^{(m)}|r^m\\
&=& \sum_{n=1}^m|b_{2n-1}|\left (\sum_{k=n}^m\left |\alpha_{2k-1}^{(2n-1)}\right |r^{2k-1} \right ). \nonumber
\eeq
%so
%\begin{equation}\label{AKP-eq5}
%\sum_{k=1}^m|a_{2k-1}|r^{2k-1} \le \sum_{n=1}^m|b_{2n-1}|\sum_{k=n}^m|\alpha_{2k-1}^{(2n-1)}|r^{2k-1}
%\end{equation}
Now for the series $\sum_{k=n}^m\left |\alpha_{2k-1}^{(2n-1)}\right |r^{2k-1}$, since $|\omega^n(z)/z^n| < 1$ for any $n\ge 1$,
 Lemma \ref{AKP-lem1} yields that
$$
\sum_{k=n}^m\left |\alpha_{2k-1}^{(2n-1)}\right |r^{2(k-n)}\le \sum_{k=n}^\infty \left |\alpha_{2k-1}^{(2n-1)}\right |r^{2(k-n)} \le 1
~\mbox{ for $\ds r\le \frac1{\sqrt3}$}.
$$
Consequently,
$$
\sum_{k=n}^m\left |\alpha_{2k-1}^{(2n-1)}\right |r^{2k-1} \le r^{2n-1}~\mbox{ for $\ds r\le \frac1{\sqrt3}$}
$$
and therefore, \eqref{AKP-eq5} reduces to
$$
\sum_{k=1}^m|a_{2k-1}|r^{2k-1} \le \sum_{n=1}^m|b_{2n-1}|r^{2n-1}~\mbox{ for $\ds r\le \frac1{\sqrt3}$, and for each $m\ge1$.}
$$
Finally, allowing $m\to\infty$ in the last inequality gives the desired inequality \eqref{AKP-eq6}.
 $\hfill \Box$

%\begin{equation}\label{AKP-eq6}
%\sum_{k=1}^\infty|a_{2k-1}|r^{2k-1} \le \sum_{n=1}^\infty|b_{2n-1}|r^{2n-1}~\mbox{ for $\ds r\le \frac1{\sqrt3}$}
%\end{equation}
%Theorem \ref{AKP-th1} is proved

%%%%%%%%%%%%%%%%%%%%%%%%%%%%%%%%%%%%%%%%%%%%%%%%%%%%%%%%%%%%%%%%%%%%%%%%%%%%%%%%%%%%%%%%%%%

%\section{ The improved version of the classical Bohr's inequality}\label{AKP-sec3}
%
%The classical Bohr inequality is not sharp for any individual function. Namely, it is easy to show that for any given function the Bohr radius is $>1/3$.
%Here is the sharp result which shows that $1/3$ cannot be improved even in individual case.
%
%\bthm\label{AKP-th2}
%Let $f(z)=\sum_{k = 0}^\infty a_kz^k$ be an analytic function in $\ID$ and $|f(z)|\le 1$ for all $z \in \ID$. Then the following
%sharp inequalities holds\begin{equation}\label{AKP-eq7}
%\frac{1-(1+|a_0|-|a_0|^2)r}{1-|a_0|r} +\sum_{k=1}^\infty |a_k|r^k \le 1 ~\mbox{ for all $r \le 1/3$},
%\end{equation}
%The function $f(z)=(z+a_0)/(1+\overline{a_0}z)$ shows that the equality holds for all $a_0 \in \ID$ and $r \leq 1/3$ .
%\ethm
%
%Proof. It is easy to see that the function $f$ is subordinated to $g=(z+a_0)/(1+\overline{a_0}z)$. Now the proof follows from Theorem B and the fact that
%$$
%\sum_{k=1}^\infty |b_k|r^k=r\frac{1-|a_0|^2}{1-r|a_0|}.
%$$

%%%%%%%%%%%%%%%%%%%%%%%%%%%%%%%%%%%%%%%%%%%%%%%%%%%%%%%%%%%%%%%%%%%%%%%%%%%%%%%%%%%%%%%%%%%

\subsection{Improved version of Bohr's inequality for  harmonic mappings}\label{AKP-sec4}

For the proof of the new version of Bohr's inequality for harmonic mappings, namely,  Theorem \ref{AKP-th3}, we need the following lemma.

\blem\label{Bohr-lem2}
Suppose that $f(z)= h(z) + \overline{g(z)} =\sum_{n = 0}^\infty a_nz^n + \overline{\sum_{n = 1}^\infty b_nz^n}$
is harmonic such that $|g'(z)| \le k|h'(z)|$ in  $\ID$ and for some $k \in [0,1]$, where $|h(z)|\le 1$  in $\ID$. Then
\elem
$$
\sum_{n=1}^\infty |a_n|r^n + \sum_{n=1}^\infty |b_n|r^n \leq (1+k)r\frac{1-|a_0|^2}{1-r|a_0|}\ \ \text{for all}\ \ r\le \frac13 .
$$
\bpf It suffices to assume that $|h(z)|< 1$ in $\ID$, and thus, hypotheses imply that the function $h$ is subordinate to $\varphi$, where
$$\varphi (z)=\frac{z+a_0}{1+\overline{a_0}z} = a_0+\sum_{k=1}^\infty \varphi_kz^k, \quad \varphi_k= (-1)^{k-1}(1-|a_0|^2)a_0^{k-1},
$$
and it is easy to see that
$$\sum_{k=1}^\infty |\varphi_k|r^k=r\frac{1-|a_0|^2}{1-r|a_0|}.
$$
Because $h(z)\prec \varphi (z)$, by using Corollary \ref{Bhowmik} and the last fact, we deduce that
\be\label{AKP-eq12}
\sum_{n=1}^\infty |a_n|r^n \leq r\frac{1-|a_0|^2}{1-r|a_0|}=1-\left [ \frac{1-r|a_0|-r( 1-|a_0|^2)}{1-r|a_0|} \right ] \ \ \text{for all}\ \ r\le \frac13.
\ee
Next, by Corollary \ref{Analog_Bhowmik}, it follows from the condition $|g'(z)| \le k|h'(z)|$ that
$$\sum_{n=1}^\infty n|b_n|r^{n-1} \le k\sum_{n=1}^\infty n|a_n|r^{n-1}.
$$
Integrating this inequality we obtain
$$\sum_{n=1}^\infty |b_n|r^{n} \le k\sum_{n=1}^\infty |a_n|r^{n}
$$
and as a consequence of it, we have
$$
\sum_{n=1}^\infty |a_n|r^n + \sum_{n=1}^\infty |b_n|r^n \leq (1+k) \sum_{n=1}^\infty |a_n|r^n \leq (1+k)r\frac{1-|a_0|^2}{1-r|a_0|},
$$
for all $0\leq r\le 1/3$, where the last inequality is a consequence of  \eqref{AKP-eq12}.
\epf

%where $B_H(r)=\sum_{n=1}^\infty |a_n|r^n + \sum_{n=1}^\infty |b_n|r^n $. % and $k=(K-1)/(K+1)$.

%New we introduce new version of Bohr's inequality for the harmonic mappings in the following theorem.
%\\

\subsection{Proof of Theorem \ref{AKP-th3}}
The proof easily follows from Lemma \ref{Bohr-lem2}.  $\hfill \Box$

%\subsection{Improved version of Bohr's inequality with constant term by function}\label{AKP-sec4}
\subsection{Proof of Theorem \ref{AKP-th5}}
By assumption $f(z)=\sum_{k=0}^\infty a_kz^k$ is analytic in $\ID$ and $|f(z)|<1$ in $\ID$.
Since $f(0)=a_0$, by assumption, the Schwarz-Pick lemma (often referred as Lindel\"{o}f's inequality) applied to the function $f$ shows that
$$|f(z)|\leq \frac{r+a}{1+ar} ~\mbox{ for $|z|=r$},
$$
where $a=|a_0|$. By using Corollary \ref{AKP-th2} we can write
$$ \sum_{k=1}^\infty |a_k|r^k \le  \frac{(1-a^2)r}{1-ar} ~ \mbox{ for all $\ds r\le \frac13$} .
$$
Combining the last two inequalities, we have
\be \label{AKP-eq13}
|f(z)|+\sum_{k=1}^\infty |a_k|r^k \le \frac{r+a}{1+ar}+ \frac{(1-a^2)r}{1-ar}  = \frac{2(1-a^2)r}{1-a^2r^2} +a
 ~ \mbox{ for all $\ds r\le \frac13$},
\ee
which is less than or equal to $1$ if
%$$
%2r\frac{1-a^2}{1-a^2r^2} +a\le 1
%$$
%so we have
\begin{equation}\label{AKP-eq9}
r^2a^2+2ra+2r-1\le 0.
\end{equation}
Solving this inequality,  we obtain that
$$ r\le r_a =\frac{\sqrt{(1+a)^2+a^2}-(1+a)}{a^2} =\frac{1}{\sqrt{(1+a)^2+a^2}\, +1+a}.
$$
We are restricted by the inequality $r(a) \leq 1/3$, which gives the condition $a\ge 2\sqrt3-3$.
This means that for $a\ge 2\sqrt3-3$ and $r\le r_a$, the desired inequality, namely \eqref{AKP-eq8}, holds.
The first part of the theorem is proved.

To show the sharpness of the radius $r_a$, we let $a=|a_0|\in[0,1)$ and consider the function
\begin{equation*}
f(z)=\frac{a_0-z}{1-\overline{a_0}z}=a_0-(1-|a_0|^2)\sum_{k=1}^\infty (\overline{a_0})^{k-1}z^k, \ \ z\in \ID.
\end{equation*}
For this function, we observe that for $z=-r$ and $a_0\geq 0$
$$
|f(z)|+\sum_{k=1}^\infty |a_k|r^k = \frac{r+a}{1+ar}+ r\frac{1-a^2}{1-ar}
$$
which shows the sharpness of $r_a$. This completes the proof of the theorem.
 $\hfill \Box$
%
%\br
%Solving this equation for $a$ gives
%$$ a\le \frac{\sqrt{2(1-r)}-1}{r} =:T(r)
%$$
%It is a simple exercise to see that
%$$
%T'(r)=-\frac{1}{r^2\sqrt{2(1-r)}}\left [ 2-r -\sqrt{2(1-r)} \right ]
%$$
%and thus, $T$ is decreasing on $[0,1]$. In particular,
%$$ T(r)\geq T(\sqrt5-2)=1 ~\mbox{ for $0<r\leq \sqrt5-2 \approx 0.2360679$}
%$$
%which establishes Remark \ref{AKP-rem1}. That is, \eqref{AKP-eq9} and hence,
%\eqref{AKP-eq8} holds for all $a<1$ whenever $0<r\leq \sqrt5-2$.
%\er

\subsection{Proof of Theorem \ref{AKP-th6}}
We follow the method of proof of Theorem \ref{AKP-th5}. Accordingly, the hypotheses imply that
$$|h(z)|\le \frac{r+a}{1+ar} , \quad |z|=r,
$$
where $a=|a_0|$, $h(0)=a_0$. The last inequality and Lemma \ref{Bohr-lem2} yield that
$$
|h(z)|+\sum_{n=1}^\infty |a_n|r^n + \sum_{n=1}^\infty |b_n|r^n  \le \frac{r+a}{1+ar} + (1+k)r\frac{1-a^2}{1-ra}
\  \text{for all}\   r\le \frac13.
$$
By making the right hand side less than or equal to $1$, we get
$$
 (1+k)r\frac{1-a^2}{1-ra}\le 1 -\frac{r+a}{1+ar} =\frac{(1-a)(1-r)}{1+ar},
$$
which upon simplification gives
\begin{equation}\label{AKP-eq10}
a(a+k+ka)r^2 + (k+2)(a+1)r -1\le 0;
\end{equation}
or equivalently,
\begin{equation}\label{AKP-eq11}
r^2(k+1)a^2+r(kr+k+2)a+r(k+2)-1\le 0.
\end{equation}
Solving the inequality (\ref{AKP-eq11}), we get that
$$
r\le r_{a,k} =\frac{B_{a,k}-(a+1)(k+2)}{2a^2(k+1)+2ak}
$$\\
where
$$
B_{a,k}= \sqrt{a^2(k^2+8k+8)+2a(k^2+6k+4)+(k+2)^2}.
$$
We have to consider those values of $a$ for which the inequality $r \leq 1/3$ holds. A little algebra shows
that the inequality $r \leq 1/3$ holds for $a\ge \alpha_k$ and hence in this case for $r\le r_{a,k}$
the desired inequality \eqref{Eq9} holds. Here $\alpha_k$ is as in the statement of Theorem \ref{AKP-th6}.

To show the sharpness of the radius $r_{a,k}$,  we consider the function
$$f(z)=h(z) +\lambda \overline{h(z)}, \quad h(z)= \frac{z+a_0}{1+\overline{a_0}z},
$$
with $\lambda \to 1$.
%$$h(z) = a_0+\sum_{k=1}^\infty \varphi_kz^k, \quad \varphi_k= (-1)^{k-1}(1-|a_0|^2)a_0^{k-1},
%$$
%with $\lambda \to 1$.
For this function, we get that (for $z=r$ and $a_0\geq 0$)
$$
|h(z)|+\sum_{n=1}^\infty |a_n|r^n + \sum_{n=1}^\infty |b_n|r^n  = \frac{r+a}{1+ar} + (\lambda+1) r\frac{1-a^2}{1-ra}
$$
and the last expression shows the sharpness of $r_{a,k}$.
This completes the proof of the theorem.
 $\hfill \Box$

%\br
%Solving the inequality (\ref{AKP-eq11}) for $a$ easily yields that
%$$ a\le \frac{\sqrt{(1-r)(k^2(1-r)+8(k+1))}+k(1-r)-2(k+1)}{2r(k+1)} =:T_k(r).
%$$
%It can be easily seen that $T_k(r)$ is a decreasing function of $r$ and so for $r\le 1/3$, we have
%$$ T_k(r)\ge T(1/3)=\frac{\sqrt{k^2+12k+12}-(2k+3)}{k+1}.
%$$
%\er

%%%%%%%%%%%%%%%%%%%%%%%%%%%%%%%%%%%%%%%%%%%%%%%%%%%%%%%%%%%%%%%%%%%

%%%%%%%%%%%%%%%%%%%%%%%%%%%%%%%%%%%%%%%%%%%%%%%%%%%%%%%%%%%%%%%%%%%%%%
\subsection*{Acknowledgments}
The authors thank the referee for his/her useful suggestions for the improvement of the presentation.
The research of I.~Kayumov  was funded by the subsidy allocated to Kazan Federal University for the state assignment in the sphere of scientific activities, project   1.13556.2019/13.1. The  work of the third author is supported by Mathematical Research Impact Centric Support of DST, India  (MTR/2017/000367).

%\subsection*{Conflict of Interest}
%The authors declare that there is no conflict of interest regarding the publication of this paper.

\end{document}